\documentclass[12pt,twoside]{article}
\usepackage{amsfonts}
\usepackage[dvips]{graphicx}
\usepackage{amsmath,fleqn}
\usepackage{amssymb}
\usepackage{xcolor}
\usepackage{mathrsfs}
\usepackage{amstext}
\usepackage{CJKutf8}
\numberwithin{figure}{section} \numberwithin{equation}{section}
\makeatletter 
\author{\small{ Lin Zhao\footnote{Corresponding author. Fax: +86-516-8912481. E-mail address: zhaolinmath@163.com}}
\medskip\\
\emph{\small{ Department of Mathematics,
China University of Mining and Technology, }}\\
\emph{\small{ Xuzhou 221116, China}}\\
\emph{\small{Corresponding author. E-mail address: zhaolinmath@163.com} }}

\begin{document}
\date{}
\title{\textbf{On some multiple solutions for a $p(x)$-Laplace equation with supercritical growth }} \maketitle
\noindent{}
\medskip\\
\indent \small{\textbf{Abstract:}  We consider the multiplicity of solutions for the $p(x)$-Laplacian problems involving the supercritical Sobolev growth via Ricceri's principle. By means of the truncation combining with De Giorgi iteration, we can extend the result about subcritical and critical growth to the supercritical growth and obtain at least three solutions for the $p(x)$ Laplacian problem.}\\
\textbf{Keywords}: Degenerate elliptic equations; a priori estimates; supercritical Sobolev exponent; truncated technique\\
\textbf{2020 Mathematics Subject Classification}: 35J70, 35B45, 35J62
\medskip\\

\section{Introduction and main result}\label{1}
~~In this paper we study the following problem
\begin{eqnarray}\label{1.1}
\left\{
  \begin{array}{ll}
    -\text{div}(|\nabla u|^{p(x)-2}\nabla u)=\lambda f(x,u)+\mu |u|^{q(x)-2}u,~~~&\text{in} ~\Omega,~~~~~~~~~~~~~~~~~~~~~~~~~~~~~~~~~ \\
    u(x)=0,~~~&\text{on} ~\partial\Omega,\\
  \end{array}
\right.
\end{eqnarray}
where $\Omega\subset \mathbb{R}^N$ is a bounded domain with smooth boundary, $p,q:\overline{\Omega}\rightarrow \mathbb{R}$ are continuous and satisfy
\begin{eqnarray}\label{1.2}
1<p_-=\min_{x\in\overline{\Omega}} p(x)\leq p_+=\max_{x\in\overline{\Omega}}p(x)<N,
\end{eqnarray}
\begin{eqnarray}
p\in L_+^\infty(\Omega)=\{p\in L^\infty(\Omega)\mid \text{ess}\inf_\Omega p(x)\geq 1\}, 
\end{eqnarray}
\begin{eqnarray}\label{1.20}
p_+<  p^\ast(x)=\frac{Np(x)}{N-p(x)}, ~q\succ p^\ast,
\end{eqnarray}
$q\succ p^\ast$ means $\inf_\Omega (q(x)-p^\ast(x))>0$, $\lambda,\mu>0$ are nonnegative parameters, and $f(x,u)$ is a Caratheodory functions satisfying
\begin{eqnarray}\label{1.3}
\sup_{s\in[-M,M]}|f(x,s)|<+\infty,~\forall M>0.
\end{eqnarray}
We suppose that the nonlinearity $f(x,u)$ satisfies the following conditions:
\begin{eqnarray}\label{1.4}
\lim_{|s|\rightarrow+\infty}\frac{f(x,s)}{|s|^{p(x)-1}}=0,~\text{uniformly~in}~ x\in \Omega;
\end{eqnarray}
\begin{eqnarray}\label{1.5}
\lim_{|s|\rightarrow0}\frac{f(x,s)}{|s|^{p(x)-1}}=0,~\text{uniformly~ in}~ x\in \Omega;
\end{eqnarray}
\begin{eqnarray}\label{1.6}
\sup_{u\in W_0^{1,p(x)}}\int_\Omega\int_0^uf(x,s)dsdx>0.
\end{eqnarray}
Let us mention here that several results have been devoted to the investigation of related problem. Fan and Deng \cite{FD} study the $p(x)$ Laplacian equation
\begin{eqnarray}\label{1.9}
\left\{
  \begin{array}{ll}
    -\text{div}(|\nabla u|^{p(x)-2}\nabla u)+a(x)|u|^{p(x)-2}u=\lambda f(x,u)+\mu g(x,u),~~~&\text{in} ~\Omega,~~~ \\
    u(x)=0,~~~&\text{on} ~\partial\Omega,\\
  \end{array}
\right.
\end{eqnarray}
with $\lambda=1$, $f$ and $g$ satisfying the subcritical growth, and obtain three solutions by the variational principle of Ricceri. Silva \cite{Si} deals with the problem (\ref{1.9}) with $a(x)=0$, $\lambda=1$, $g(x,u)=|u|^{q(x)-2}u$, $q(x)\leq p^\ast(x)$ and $f(x,u)$ satisfying the symmetry condition; the authors obtains that the problem (\ref{1.9}) has at least $m$ pair of solutions for given $m\in \mathbb{N}$ when $\mu>0$ small by the symmetric mountain pass lemma and the concentration compactness lemma of Lions \cite{L2} for variable exponent. Fan and Zhao \cite{FZY} study the existence of nodal radial solutions for the problem (\ref{1.9}) with $\lambda=1$, $\mu=0$ and $f(x,u)=|u|^{q(x)-2}u$; get a pair of solutions which have exactly $k$ nodes when $a(x)$, $p(x)$ and $q(x)$ satisfy suitable conditions. By the truncation and Moser iteration \cite{Mo}, Zhao and Zhao \cite{ZZ} study the equation (\ref{1.1}) with the supercritical growth as $p(x)$ is a constant exponent, and obtain at least three solutions. In this paper we study the problem (\ref{1.1}) with the supercritical growth by the truncated technique and De Giogri iteration \cite{FZ,FS,De}. If $a(x)=0$, $g(x,u)=|u|^{q(x)-2}u$ and $q\succ p^\ast$, the problem (\ref{1.9}) is (\ref{1.1}). For the case $q\succ p^\ast$, we can not use directly the variational techniques because the corresponding functional is not well-defined on the Sobolev space.\\

Obviously, The condition (\ref{1.6}) implies
\begin{eqnarray}\label{1.10.0}
\theta:=\sup_{\Phi(u)\neq 0}\frac{J(u)}{\Phi(u)}>0,
\end{eqnarray}
where $\Phi(u)=\int_\Omega \frac{1}{p(x)}|\nabla u|^{p(x)}dx$ and $J(u)=\int_\Omega\int_0^uf(x,s)dsdx$.

The main result in this paper is stated as follows:\\
\textbf{Theorem 1.1.} Let (\ref{1.2})--(\ref{1.6}) hold, and $\theta$ be given by (\ref{1.10.0}). Then for each compact
interval $[a,b]\subset(\frac{1}{\theta},+\infty)$ there exists $\gamma>0$ with the following property: for every $\lambda\in[a,b]$,
there exists $\delta>0$ such that, for each $\mu\in[0,\delta]$, the problem (\ref{1.1}) has at least three solutions in $W_0^{1,p(x)}(\Omega)\bigcap L^{\infty}(\Omega)$, whose
$W_0^{1,p(x)}(\Omega)$-norms are less than $\gamma$.\\
\textbf{Remark 1.2.} For the condition $q\succ p^\ast$ in problem (\ref{1.1}), this case is not the natural growth condition on $u$. One can only get the Sobolev embedding $W_0^{1,p(x)}(\Omega)\hookrightarrow L^{q(x)}(\Omega)$, $q(x)\leq p^\ast(x)=\frac{Np(x)}{N-p(x)}$. In this paper we are interested in existence and multiplicity of solutions of the problem (\ref{1.1}) with $q\succ  p^\ast$. One can not directly use variational methods, as there does not exist the embedding $W_0^{1,p(x)}(\Omega)\hookrightarrow L^{q(x)}(\Omega)$, $q\succ p^\ast$. There is a difficulty for the $p(x)$ Laplacian equation with supercritical growth. The difficulty is the lack of the Sobolev embedding and compactness, and we overcome this difficulty by the truncated technique and De Giogri iteration \cite{FZ,FS,De}.

\section{Preliminaries}\label{sec2}
Let $\Omega$ be an open domain of $\mathbb{R}^N$, denote
\begin{equation}
L_+^\infty(\Omega)=\{p\in L^\infty(\Omega)\mid \text{ess}\inf_\Omega p(x)\geq 1\}.
\end{equation}
For $p\in L_+^\infty(\Omega)$, define the space
\begin{equation}
L^{p(x)}(\Omega)=\{u:\Omega\rightarrow \mathbb{R}~\text{measurable}\mid \int_\Omega |u(x)|^{p(x)}dx<\infty\}
\end{equation}
with the norm
\begin{equation}
\|u\|_{L^{p(x)}}=\inf\{t>0\mid \int_\Omega|\frac{u(x)}{t}|^{p(x)}dx\leq 1 \},
\end{equation}
and Sobolev space with variable exponent
\begin{equation}
W^{1,p(x)}(\Omega)=\{u\in L^{p(x)}(\Omega)\mid  |\nabla u|\in L^{p(x)}\}
\end{equation}
with the norm
\begin{equation}
\|u\|_{W^{1,p(x)}}=\|u\|_{L^{p(x)}}+\|\nabla u\|_{L^{p(x)}}.
\end{equation}
We denote $W_0^{1,p(x)}(\Omega)=\overline{C_0^\infty(\Omega)}^{\|\cdot\|_{}}$, with the norm $\|u\|=\|\nabla u\|_{L^{p(x)}}$.\\
\textbf{Lemma 2.1} (\cite{FD,Si,FZ2}). Let $\rho(u)=\int_\Omega |u(x)|^{p(x)}dx$. For $u,u_k\in L^{p(x)}(\Omega)$, we have
\begin{eqnarray}
&(1)&~\text{For}~u\neq0,\|u\|_{L^{p(x)}}=t\Leftrightarrow \rho(\frac{u}{t})=1. \nonumber\\
&(2)&~\|u\|_{L^{p(x)}}<1(=1;>1)\Leftrightarrow \rho(u)<1(=1;>1).\nonumber\\
&(3)&~\|u\|_{L^{p(x)}}>1\Rightarrow \|u\|_{L^{p(x)}}^{p_-}\leq \rho(u)\leq\|u\|_{L^{p(x)}}^{p_+}.\nonumber\\
&(4)&~\|u\|_{L^{p(x)}}<1\Rightarrow \|u\|_{L^{p(x)}}^{p_+}\leq \rho(u)\leq\|u\|_{L^{p(x)}}^{p_-}.\nonumber\\
&(5)&~\|u_k\|_{L^{p(x)}}=o_k(1)\Leftrightarrow \rho(u_k)=o_k(1).\nonumber\\
&(6)&~\|u_k\|^{-1}_{L^{p(x)}}=o_k(1)\Leftrightarrow (\rho(u_k))^{-1}=o_k(1).\nonumber\\
\end{eqnarray}
\textbf{Lemma 2.2}(\cite{FZ2}). The spaces $L^{p(x)}(\Omega)$ and $W^{1,p(x)}(\Omega)$ are separable Banach spaces, and they are reflexive as $p_->1$.\\
\textbf{Lemma 2.3}(\cite{FD,FZ2}). Let $\Omega$ be an open bounded domain in $\mathbb{R}^N$ with the cone property, $p$ satisfy (\ref{1.2}), $q\in C^0(\mathbb{R}^N)$ and $p^\ast\succ q$, then the embedding $W^{1,p(x)}(\Omega)\hookrightarrow L^{q(x)}(\Omega)$ is compact, where  $p^\ast \succ q$ means $\inf_\Omega (p^\ast(x)-q(x))>0$.\\
\textbf{Lemma 2.4} (\cite{FS}). Suppose a sequence $a_i$, $i=0,1,...$ of nonnegative numbers satisfies the recursion relation
\begin{equation}
a_{i+1}\leq c b^i a_i^{1+\eta},~i=0,1,2...,
\end{equation}
where $c$, $b$ and $\eta$ are positive constants and $b>1$. If $a_0\leq c^{-\frac{1}{\eta}} b^{-\frac{1}{\eta^2}}$, then $a_i\rightarrow 0$ as $i\rightarrow\infty$.\\

If $X$ is a real Banach space, we can denote by $\mathfrak{W}_X$ the class of all functionals $\phi: X\rightarrow\mathbb{R}$ possessing the following property: if $\{u_n\}\subset X$ is a sequence converging weakly to $u\in X$ and $\liminf_{n\rightarrow\infty}\phi(u_n)\leq\phi(u)$, then $\{u_n\}$ has a subsequence converging strongly to $u$. For instance, if $X$ is an uniformly convex Banach space and $g:[0,+\infty)\rightarrow\mathbb{R}$ is a continuous, strictly increasing function, then the functional $u\rightarrow g(\|u\|)$ belongs to the class $\mathfrak{W}_X$.\\
\textbf{Lemma 2.5} (\cite{R1,R2}). Let $X$ be a separable and reflexive real Banach space; 
$\Phi:X\rightarrow\mathbb{R}$ a sequentially weakly lower semicontinuous $C^{1}$ functional, belonging to
$\mathfrak{W}_X$, bounded on each bounded subset of $X$ and whose derivative admits a continuous inverse on $X^{\ast}$;
$J:X\rightarrow \mathbb{R}$ a $C^1$ functional with compact derivative. Assume that $\Phi$ has a strict local minimum $x_0$ with $\Phi(x_0)=J(x_0)=0$. Finally, setting
$$\alpha=\max\{0,\limsup_{\|x\|\rightarrow\infty}\frac{J(x)}{\Phi(x)}, \limsup_{x\rightarrow x_0}\frac{J(x)}{\Phi(x)}\},$$
$$\beta=\sup_{x\in \Phi^{-1}(0,+\infty)}\frac{J(x)}{\Phi(x)},$$
assume that $\alpha <\beta$. Then, for each compact interval $[a,b]\subseteq I=(\frac{1}{\beta},\frac{1}{\alpha})$ (with the conventions $\frac{1}{0}=+\infty$, $\frac{1}{+\infty}=0$), there exists $\rho>0$ with the following
property: for every $\lambda\in [a,b]$ and every $C^1$ functional $\Psi: X\rightarrow\mathbb{R}$ with compact derivative,
there exists $\delta_0>0$ such that, for each $\mu\in [0,\delta_0]$ the equation
\begin{equation*}
\Phi'(u)=\lambda J'(u)+\mu\Psi'(u)
\end{equation*}
has at least three solutions in $X$ whose norms are less than $\gamma$, i.e. $\|u_i\|\leq\gamma$, $i=1,2,3$.\\
\textbf{Remark 2.6.} From the proof of Lemma 2.5, Ricceri \cite{R2} shows that $\rho$ does not depend on $\mu$. One can refer to \cite{R2} for detail.\\
\textbf{Lemma 2.7.} Let $f$ satisfy (\ref{1.3}) and (\ref{1.4}). Then for every $\lambda\in(0,+\infty)$, the functional $\Phi-\lambda J$ is  sequentially weakly
lower continuous and coercive on $W_0^{1,p(x)}(\Omega)$, and has a global minimizer $v_\lambda$, where $\Phi(u)= \int_{\Omega}\frac{1}{p(x)}|\nabla u|^{p(x)}dx$ and $J(u)= \int_{\Omega}\int_0^uf(x,s)dsdx$.    \\
\textbf{Proof.} Let us fix $\lambda\in (0,+\infty)$. By
(\ref{1.4}), for $\forall\varepsilon>0$, there exists $M_0>0$, such that
\begin{equation*}
|f(x,s)|\leq\varepsilon|s|^{p(x)-1}~~~\text{as}~|s|\geq M_0.
\end{equation*}
Considering this inequality and (\ref{1.3}), we may find a constant $C_0>0$, such that
\begin{equation}
|f(x,s)|\leq C_0+\varepsilon|s|^{p(x)-1}~~~\text{as}~s\in\mathbb{R},
\end{equation}
which implies
\begin{equation}\label{2.9}
|F(x,u)|\leq C_0|u|+\frac{\varepsilon}{p(x)}|u|^{p(x)}.
\end{equation}
Thus, for $u\in W_{0}^{1,p(x)}(\Omega)$, we obtain
\begin{eqnarray}
\Phi(u)-\lambda J(u)&=&\int_{\Omega}\frac{1}{p(x)}|\nabla u|^{p(x)}dx-\lambda\int_{\Omega}F(x,u)dx\nonumber\\
&\geq&\int_{\Omega}\frac{1}{p(x)}|\nabla u|^{p(x)}dx-\lambda\int_{\Omega}(\frac{\varepsilon}{p(x)}|u|^{p(x)}+C_0|u|)dx\nonumber\\
&\geq&\frac{1}{p_+}\parallel u\parallel^{p(x)}-\frac{\lambda\varepsilon}{p_-}\int_{\Omega}|u|^{p(x)}dx-\lambda C_0\int_{\Omega}|u|dx\nonumber\\
&\geq&\frac{1}{p_+}\parallel u\parallel^{p(x)}-\frac{\lambda\varepsilon C_p^{p}}{p_-}\parallel u\parallel^{p(x)}-\lambda C_0C_1\parallel u\parallel\nonumber,
\end{eqnarray}
where $\|u\|_{L^{p(x)}}\leq C_p\|u\|$ and $\|u\|_{L^1}\leq C_1\|u\|$, with constants $C_1,C_p>0$.
Since $p(x)>1$, $\varepsilon>0$ is small enough, so we have $\Phi(u)-\lambda J(u)\rightarrow+\infty$ as $\parallel u\parallel\rightarrow\infty$. Hence $\Phi-\lambda J$ is coercive.\\
Since the embedding $W_0^{1,p(x)}(\Omega)\hookrightarrow L^{p(x)}(\Omega)$ is compact and (\ref{2.9}), $J$ is weakly continuous. Obviously, $\Phi(u)=\int_{\Omega}\frac{1}{p(x)}|\nabla u|^{p(x)}dx\in\mathfrak{W}_X$ is weakly lower semicontinuous on $W_0^{1,p(x)}(\Omega)$. We can deduce that $\Phi-\lambda J$ is sequentially weakly lower semicontinuous. So $\Phi-\lambda J$ has a global minimizer $v_\lambda$.\\
The proof is complete.\\
Next, we will show that $\Phi-\lambda J$ has a strictly local, not global minimizer for some $\lambda$, when $f$ satisfies (\ref{1.3})--(\ref{1.6}).\\
\textbf{Lemma 2.8.} Let $f$ satisfy (\ref{1.3})--(\ref{1.6}). Then\\
$(i)$ $0$ is a strict local minimizer of the functional $\Phi-\lambda J$ for $\lambda\in(0,+\infty)$.\\
$(ii)$ $0\neq v_\lambda$, i.e., $0$ is not the global minimizer $v_\lambda$ for $\lambda\in (\frac{1}{\theta},+\infty)$, where $v_\lambda$ is given by Lemma 2.7.\\
\textbf{Proof.}
Firstly, we prove
\begin{equation}
\lim_{\|u\|\rightarrow0}\frac{J(u)}{\Phi(u)}=0.
\end{equation}
In fact, by (\ref{1.5}), for $\forall\varepsilon>0$, $\exists \delta>0$, such that
\begin{equation}\label{2.11}
|f(x,u)|\leq\varepsilon|u|^{p(x)-1},~~\text{as}~|u|<\delta.
\end{equation}
Considering the inequality (\ref{2.11}) and (\ref{1.4}), we have
\begin{equation}
|f(x,u)|\leq\varepsilon|u|^{p(x)}+|u|^{r(x)},~~\text{for}~x\in\Omega,~u\in\mathbb{R},
\end{equation}
with fixed $r(x)\in(p(x),p^\ast(x))$. By continuous embedding, we have
\begin{eqnarray}
|J(u)|&\leq& \int_\Omega |F(x,u)|dx\nonumber\\
&\leq& \frac{\varepsilon}{p_-}\int_\Omega |u|^{p(x)}dx +\frac{1}{r_-}\int_\Omega |u|^{r(x)}dx\nonumber\\
&\leq&\frac{\varepsilon C_p}{p_-} \|u\|^{p(x)}+\frac{C_r}{r_-}\|u\|^{r(x)}.
\end{eqnarray}
This implies
\begin{equation}
\lim_{\|u\|\rightarrow0}\frac{J(u)}{\Phi(u)}=0.
\end{equation}
Next, we will prove $(i)$ and $(ii)$.\\
$(i)$ For $\lambda\in(0,+\infty)$, since $\lim_{\|u\|\rightarrow0}\frac{J(u)}{\Phi(u)}=0<\frac{1}{\lambda}$ and $\Phi(u)>0$ for each $u\neq0$ in some neighborhood $U$ of $0$, there exists a neighborhood $V\subseteq U$ of $0$ such that $\Phi(u)-\lambda J(u)>0$ for all $u\in V\setminus\{0\}$. Hence, $0$ is a strict local minimum of $\Phi-\lambda J$.\\
$(ii)$ For $\lambda\in(\frac{1}{\theta},+\infty)$, from the definition of $\theta$ (\ref{1.10.0}), there exists $u^\ast\in W_0^{1,p(x)}(\Omega)$, with $\min\{\Phi(u^\ast)$, $J(u^\ast)\}>0$, such that $\frac{J(u^\ast)}{\Phi(u^\ast)}>\frac{1}{\lambda}$, i.e. $\Phi(u^\ast)-\lambda J(u^\ast)<0=\Phi(0)-\lambda J(0)$. So $0\neq v_\lambda$ is not a global minimum of $\Phi-\lambda J$.\\
The proof is complete.

\section{Proof of Theorem 1.1}\label{sec3}
Our main difficulty to prove the existence solutions of the problem (\ref{1.1}) by using variational methods is that the associated to functional is not well defined on $W^{1,p(x)}_0(\Omega)$ for $q\succ  p^\ast$. We follow idea from \cite{ZZ,FZ}. Let the truncation of $|u|^{q(x)-2}u$ be given by
\begin{eqnarray}\label{3.1.0}
g_K(x,u)=
\left\{
  \begin{array}{ll}
    |u|^{q(x)-2}u,~~&\text{if}~0\leq|u|\leq K,~~~~~~~~~~ \hbox{} \\
    K^{q(x)-p_+}|u|^{p_+-2}u,~~~~~~~~~~~~~~~~~&\text{if}~|u|\geq K, \hbox{}
  \end{array}
\right.
\end{eqnarray}
where $K\geq2$ is a real number, whose value will be fixed later. Then $g_K(x,u)$ satisfies the subcritical growth
\begin{eqnarray}\label{3.2}
|g_K(x,u)|\leq K^{q(x)-p_+}|u|^{p_+-1}.
\end{eqnarray}
Setting $h_K(x,u)=\lambda f(x,u)+\mu g_K(x,u)$, we study the truncated problem associated to $h_K$
\begin{eqnarray}\label{3.3}
\left\{
  \begin{array}{ll}
    -\text{div}(|\nabla u|^{p(x)-2}\nabla u)=h_K(x,u)=\lambda f(x,u)+\mu g_K(x,u),~~&\text{in} ~\Omega,~~\\
    u(x)=0,&\text{on} ~\partial\Omega.\\
  \end{array}
\right.
\end{eqnarray}
\textbf{Definition 3.1.} We say that $u\in W_{0}^{1,p(x)}(\Omega)$ is a weak
solution of the truncated problem (\ref{3.3}) if
\begin{eqnarray}\label{3.4}
\int_{\Omega}|\nabla u|^{p(x)-2}\nabla u\cdot\nabla \varphi
dx=\int_{\Omega}h_K(x,u)\varphi
dx
\end{eqnarray}
for every $\varphi\in W_{0}^{1,p(x)}(\Omega)$. We introduce
$$\Psi (u)=\int_\Omega\int_0^ug_K(x,s)dsdx.$$
It is easy to see that the critical points of the functional $\Phi-\lambda J-\mu \Psi$ are the weak solutions of the problem (\ref{3.4}), where $\Phi(u)= \int_{\Omega}\frac{1}{p(x)}|\nabla u|^{p(x)}dx$ and $J(u)= \int_{\Omega}\int_0^uf(x,s)dsdx$. From (\ref{1.5}), for small $\varepsilon>0$, $\exists \delta >0$ such that the function $f$ satisfies
\begin{eqnarray}
|f(x,s)|\leq \varepsilon |s|^{p(x)-1},
\end{eqnarray}
for $|s|<\delta$. For each compact interval $[a,b]\subset(\frac{1}{\theta},+\infty)$, $\lambda\in [a,b]$, considering (\ref{1.4}) and $|g_K(x,u)|\leq K^{q(x)-p_+}|u|^{p_+-1}$, we can choose constant $C>0$ such that
\begin{eqnarray}\label{3.6}
|h_K(x,s)|\leq C|s|^{p(x)-1}+\mu K^{q(x)-p_+}|s|^{p_+-1},
\end{eqnarray}
for $\forall s\in \mathbb{R}$. Hence, $h_K(x,u)$ and $g_K(u)$  satisfy the subcritical growth. $\Psi(u)$ has a compact derivative in $W^{1,p(x)}_0(\Omega)$. By Lemma 2.2, Lemma 2.3, Lemma 2.7 and Lemma 2.8, all the hypotheses
of Lemma 2.5 are satisfied. Then there exists $\gamma>0$, with the
following the property: for every $\lambda\in[a,b]$, there exists
$\delta_0>0$, such that, for $\mu\in[0,\delta_0]$, the problem
(\ref{3.3}) has at least three solutions $u_0$, $u_1$ and $u_2$ in
$W_0^{1,p(x)}(\Omega)$ whose $W_0^{1,p(x)}(\Omega)$-norms are less than
$\gamma$, i.e., $\|u_i\|_{W_0^{1,p(x)}}\leq\gamma$, $i=0,1,2$, where $\gamma$ is depend on $\lambda$, but not depend on $\mu$ (see (5) of \cite{R2} for details). If the
three solutions $u_i$, $i=0,1,2$, satisfy
\begin{equation}\label{3.7}
|u_i(x)|\leq K, ~\text{a.e.}~x\in\Omega, ~i=0,1,2,
\end{equation}
then in view of the definition $g_K$ (\ref{3.1.0}), we have $h_K(x,u)=\lambda f(x,u)+\mu|u|^{q(x)-2}u$ and therefore $u_i$, $i=0,1,2$, are also solutions of the original problem (1.1). Thus, in order to prove Theorem 1.1, it suffices to show that exists $\delta>0$, such that for $\mu\in[0,\delta]$, the solutions obtained by Lemma 2.5 satisfy the inequality (\ref{3.7}).\\
\\
\textbf{Step 1.} We will firstly prove that $\forall$ any ball $B_R(x_0)\subset \subset\Omega$, $R\leq1$, $\forall$ real number $l\geq1$ and $\forall$ $t,s$ satisfying $0\leq t<s\leq1$ the solutions $u=u_i$, $i=0,1,2$ satisfy the following inequality
\begin{eqnarray}\label{3.8}
\int_{A_{l,t}}|\nabla u|^{p_-}dx\leq (C'+C'' \mu K^{q_+-p_+})(\int_{A_{l,s}}|\frac{u-l}{s-t}|^{p_-^\ast}dx+(l^{p_+}+1)|A_{l,s}|),
\end{eqnarray}
where $p_-^\ast=\frac{Np_-}{N-p_-}$, $C',C''>0$ are constants and $A_{l,s}=\{x\in B_s\mid u(x)>l\}$. \\
For arbitrary balls $\overline{B}_t(x')\subset B_s(x') \subset B_R(x_0)$, some $x'\in \Omega$, let $\xi$ be a $C^\infty$ function such that $0\leq \xi \leq 1$, $\text{supp}\xi\subset B_s$, $\xi=1$ on $B_t$, $|\nabla \xi|\leq \frac{2}{s-t}$. For $l\geq 1$ taking $\varphi=\xi^{p_+}\max\{u-l,0\}\in W^{1,p(x)}_0(\Omega)$ as a test function in the problem (\ref{3.4}), we get
\begin{eqnarray}\label{3.9}
\int_{\Omega}|\nabla u|^{p(x)-2}\nabla u\nabla\varphi dx=\int_{\Omega}[\lambda f(x,u)\varphi +\mu g_K(x,u)\varphi ]dx,
\end{eqnarray}
and
\begin{eqnarray}\label{3.10}
\int_{\Omega}|\nabla u|^{p(x)-2}\nabla u\nabla\varphi dx=\int_{A_{l,s}}|\nabla u|^{p(x)}\xi^{p_+}dx +p_+\int_{A_{l,s}}|\nabla u|^{p(x)-2}\xi^{p_+-1}(u-l)\nabla u \nabla\xi dx.\nonumber\\
\end{eqnarray}
From (\ref{3.6}), we have
\begin{eqnarray}\label{3.11}
&&\int_{\Omega}[\lambda f(x,u)\varphi +\mu g_K(x,u)\varphi ]dx\nonumber\\
&=&\int_{A_{l,s}}[\lambda f(x,u)\varphi +\mu g_K(x,u)\varphi ]dx\nonumber\\
&\leq&\lambda C\int_{A_{l,s}}|u|^{p(x)-1}\xi^{p_+}(u-l)dx+\mu \int_{A_{l,s}} K^{q(x)-p_+}|u|^{p_+-1}\xi^{p_+}(u-l)dx.
\end{eqnarray}
By (\ref{3.9})--(\ref{3.11}), it follows that
\begin{eqnarray}\label{3.12}
\int_{A_{l,s}}|\nabla u|^{p(x)}\xi^{p_+}dx&\leq&-p_+\int_{A_{l,s}}|\nabla u|^{p(x)-2}\xi^{p_+-1}(u-l)\nabla u \nabla\xi dx\nonumber\\
&&+\lambda C\int_{A_{l,s}}|u|^{p(x)-1}\xi^{p_+}(u-l)dx \nonumber \\
&&+\mu \int_{A_{l,s}} K^{q(x)-p_+}|u|^{p_+-1}\xi^{p_+}(u-l)dx \nonumber  \\
&\leq&p_+\int_{A_{l,s}}|\nabla u|^{p(x)-1}\xi^{p_+-1}(u-l)| \nabla\xi| dx  \nonumber   \\
&&+\lambda C\int_{A_{l,s}}|u|^{p(x)-1}\xi^{p_+}(u-l)dx   \nonumber  \\
&&+\mu \int_{A_{l,s}} K^{q(x)-p_+}|u|^{p_+-1}\xi^{p_+}(u-l)dx.
\end{eqnarray}
Let us estimate the terms of the right hand side of (\ref{3.12}),
\begin{eqnarray}\label{3.13}
&&p_+\int_{A_{l,s}}|\nabla u|^{p(x)-1}\xi^{p_+-1}(u-l) |\nabla\xi| dx  \nonumber   \\
&\leq& p_+[\int_{A_{l,s}} \frac{p(x)-1}{p(x)}\varepsilon^\frac{p(x)}{p(x)-1}|\nabla u|^{p(x)}\xi^\frac{(p_+-1)p(x)}{p(x)-1}dx \nonumber    \\
&&+\int_{A_{l,s}}\frac{1}{p(x)}\varepsilon^{-p(x)}|\nabla\xi|^{p(x)}|u-l|^{p(x)} dx]  \nonumber   \\
&\leq&(p_+-1)\varepsilon^{\frac{p_+}{p_+-1}}\int_{A_{l,s}}|\nabla u|^{p(x)}\xi^{p_+}dx+\frac{p_+}{p_-}\varepsilon^{-p_-}2^{p_+}\int_{A_{l,s}}|\frac{u-l}{s-t}|^{p(x)}dx  \nonumber  \\
&\leq&\frac{1}{2}\int_{A_{l,s}}|\nabla u|^{p(x)}\xi^{p_+}dx+C_{p,\varepsilon}(\int_{A_{l,s}}\frac{p(x)}{p_-^\ast}|\frac{u-l}{s-t}|^{p_-^\ast}dx +\frac{p_-^\ast-p(x)}{p_-^\ast}|A_{l,s}|)  \nonumber  \\
&\leq&\frac{1}{2}\int_{A_{l,s}}|\nabla u|^{p(x)}\xi^{p_+}dx+C_{p,\varepsilon}(\int_{A_{l,s}}|\frac{u-l}{s-t}|^{p_-^\ast}dx +|A_{l,s}|),
\end{eqnarray}
where we use the Young's inequality with $\varepsilon$ and take $\varepsilon\in(0,1)$ such that $(p_+-1)
\varepsilon^\frac{p_+}{p_+-1}=\frac{1}{2}$ and denote $C_{p,\varepsilon}=\frac{p_+2^{p_+}}{p_-}(2p_+-2)^\frac{(p_+-1)p_-}{p_+}$. Using Young's inequality, we get
\begin{eqnarray}\label{3.14}
&&\lambda C\int_{A_{l,s}}|u|^{p(x)-1}\xi^{p_+}|u-l|dx \nonumber\\
&\leq&\lambda C[\int_{A_{l,s}}\frac{1}{p(x)}|u-l|^{p(x)}dx+\int_{A_{l,s}}\frac{p(x)-1}{p(x)}|u|^{p(x)}dx ]\nonumber\\
&\leq&\frac{\lambda C}{p_+}\int_{A_{l,s}}|u-l|^{p(x)}dx+\frac{\lambda C(p_+-1)2^{p_+-1}}{p_+}[\int_{A_{l,s}}|u-l|^{p(x)}dx+l^{p_+}|A_{l,s}|]\nonumber\\
&=&\frac{\lambda C}{p_+}[((p_+-1)2^{p_+-1}+1)\int_{A_{l,s}}|u-l|^{p(x)}dx+(p_+-1)2^{p_+-1}l^{p_+} |A_{l,s}|]\nonumber\\
&\leq&\frac{\lambda C}{p_+}[2^{p_+-1}p_+(\int_{A_{l,s}}\frac{p(x)}{p_-^\ast}|\frac{u-l}{s-t}|^{p_-^\ast}dx +\int_{A_{l,s}}\frac{p_-^\ast-p(x)}{p_-^\ast}dx) +2^{p_+-1}p_+l^{p_+} |A_{l,s}|]\nonumber\\
&\leq& \lambda C 2^{p_+-1}[ \int_{A_{l,s}}|\frac{u-l}{s-t}|^{p_-^\ast}dx+(l^{p_+}+1) |A_{l,s}|].
\end{eqnarray}
Considering $K\geq 2$ and $0\leq\xi\leq1$, we get
\begin{eqnarray}\label{3.15}
&&\mu\int_{A_{l,s}}K^{q(x)-p_+}|u|^{p_+-1}\xi^{p_+}|u-l|dx\nonumber\\
&\leq&\mu K^{q_+-p_+}\int_{A_{l,s}}|u|^{p_+-1}|u-l|dx\nonumber\\
&\leq&\mu K^{q_+-p_+}(\int_{A_{l,s}}|u|^{p_+}dx+\int_{A_{l,s}}|u-l|^{p_+}dx)\nonumber\\
&\leq&\mu K^{q_+-p_+}[2^{p_+-1}(\int_{A_{l,s}}|u-l|^{p_+}dx+l^{p_+}|A_{l,s}|)+\int_{A_{l,s}}|u-l|^{p_+}dx]\nonumber\\
&\leq&\mu K^{q_+-p_+}2^{p_+}(\int_{A_{l,s}}|\frac{u-l}{s-t}|^{p_+}dx+l^{p_+}|A_{l,s}|)\nonumber\\
&\leq&\mu K^{q_+-p_+}2^{p_+}(\int_{A_{l,s}}\frac{p_+}{p_-^\ast}|\frac{u-l}{s-t}|^{p_-^\ast}dx +\frac{p_-^\ast-p_+}{p_-^\ast}|A_{l,s}| +   l^{p_+}|A_{l,s}|)\nonumber\\
&\leq&\mu K^{q_+-p_+}2^{p_+}(\int_{A_{l,s}}|\frac{u-l}{s-t}|^{p_-^\ast}dx + (l^{p_+}+1)|A_{l,s}|),
\end{eqnarray}
where $0<t< s \leq 1$. From (\ref{3.12})--(\ref{3.15}), it follows that
\begin{eqnarray}\label{3.16.0}
&&\int_{A_{l,t}}|\nabla u|^{p_-}dx\nonumber\\
&\leq& \int_{A_{l,t}}\frac{p_-}{p(x)}|\nabla u|^{p(x)}dx+\frac{p(x)-p_-}{p(x)}|A_{l,t}|\nonumber\\
&\leq&\int_{A_{l,t}} |\nabla u|^{p(x)}dx+|A_{l,t}|\nonumber\\
&\leq&\int_{A_{l,s}} |\nabla u|^{p(x)}\xi^{p_+}dx+|A_{l,s}|\nonumber\\
&\leq& 2C_{p,\varepsilon}(\int_{A_{l,s}}|\frac{u-l}{s-t}|^{p_-^\ast}dx+|A_{l,s}|)\nonumber\\
&&+(C\lambda 2^{p_+}+\mu K^{q_+-p_+}2^{p_++1})[\int_{A_{l,s}}|\frac{u-l}{s-t}|^{p_-^\ast}dx+(l^{p_+}+1)|A_{l,s}|]+|A_{l,s}|\nonumber\\
&\leq&(2C_{p,\varepsilon}+C\lambda 2^{p_+}+\mu K^{q_+-p_+}2^{p_++1}+1)[\int_{A_{l,s}}|\frac{u-l}{s-t}|^{p_-^\ast}dx+(l^{p_+}+1)|A_{l,s}|],
\end{eqnarray}
where $C_{p,\varepsilon}=\frac{p_+}{p_-}2^{p_+}(2p_+-2)^\frac{(p_+-1)p_-}{p_+}$. If we denote $C'=2C_{p,\varepsilon}+C\lambda 2^{p_+}+1$ and $C''=2^{p_++1}$ in (\ref{3.16.0}), we can get the inequality (\ref{3.8}).\\
\\
\textbf{Step 2.} We will show that $\exists \delta>0$ such that  $|u|\leq K$ a.e. $x\in \Omega$ for $\mu\in[0,\delta]$.\\
Let us fix $B_{R}\subset\subset\Omega$, $K\geq 2$ in \textbf{Step 1} and set
\begin{eqnarray}
\rho_i=\frac{R}{2}+\frac{R}{2^{i+1}},
\end{eqnarray}
\begin{eqnarray}
\bar{\rho}_i=\frac{\rho_i+\rho_{i+1}}{2},
\end{eqnarray}
and
\begin{eqnarray}
K_i=K(1-\frac{1}{2^{i+1}}),
\end{eqnarray}
with $i=0,1,2...$ Obviously,
\begin{eqnarray}
\rho_{i+1}\leq \bar{\rho}_i\leq \rho_i,
\end{eqnarray}
\begin{eqnarray}\label{3.21.0}
\rho_0\geq \rho_1\geq...\geq \rho_i...\rightarrow \frac{R}{2}
\end{eqnarray}
and
\begin{eqnarray}\label{3.22.0}
1\leq K_0=\frac{K}{2}\leq K_1\leq...\leq K_i...\rightarrow K
\end{eqnarray}
as $i\rightarrow\infty$.\\
We define the sequence
\begin{eqnarray}
a_i=\int_{A_{K_i,\rho_i}}|u-K_i|^{p_-^\ast}dx
\end{eqnarray}
with $p_-^\ast=\frac{N p_-}{N-p_-}$ and fix a function $\zeta(\cdot)\in C^1[0,+\infty)$, $0\leq \zeta(t)\leq 1$
\begin{eqnarray}
\zeta(t)=
\left\{
  \begin{array}{ll}
    1,~&\text{if}~t\leq\frac{1}{2},~~~~~~~~~~ \hbox{} \\
    0,~~&\text{if}~t\geq \frac{3}{4}, \hbox{}
  \end{array}\nonumber
\right.
\end{eqnarray}
and $|\zeta'(t)|\leq 1$. Setting
\begin{eqnarray}
\zeta_i(x)=\zeta (\frac{2^{i+1}}{R}(|x|-\frac{R}{2})),
\end{eqnarray}
we have
\begin{eqnarray}
\zeta_i(x)=
\left\{
  \begin{array}{ll}
    1,~&\text{}~x\in B_{\rho_{i+1}},~~~~~~~~~~ \hbox{} \\
    0,~&\text{}~x\notin B_{\bar{\rho}_i}.\hbox{}
  \end{array}\nonumber
\right.
\end{eqnarray}
Next, we will show that $a_i$ satisfies the inequality $a_{i+1}\leq c b^i a_i^{1+\eta}$.
\begin{eqnarray}\label{3.25}
a_{i+1}&=& \int_{A_{K_{i+1},\rho_{i+1}}}|u-K_{i+1}|^{p_-^\ast}dx\nonumber\\
&\leq& \int_{A_{K_{i+1},\bar{\rho}_{i}}}|(u-K_{i+1})\zeta_i|^{p_-^\ast}dx\nonumber\\
&\leq&\int_{B_R}|(u-K_{i+1})\zeta_i|^{p_-^\ast}dx\nonumber\\
&\leq& S^{-\frac{p_-^\ast}{p_-}}[\int_{B_R}|\nabla ((u-K_{i+1})\zeta_i)|^{p_-}dx]^{\frac{p_-^\ast}{p_-}} \nonumber\\
&\leq& S^{-\frac{p_-^\ast}{p_-}}[\int_{B_R}|\zeta_i\nabla u+(u-K_{i+1})\nabla\zeta_i|^{p_-}dx]^{\frac{p_-^\ast}{p_-}} \nonumber\\
&\leq& S^{-\frac{p_-^\ast}{p_-}} 2^{p_-^\ast}  [\int_{A_{K_{i+1}},\bar{\rho}_i}|\zeta_i\nabla u|^{p_-}dx+ \int_{A_{K_{i+1}},\rho_i}|(u-K_{i+1})\nabla\zeta_i|^{p_-}dx]^{\frac{p_-^\ast}{p_-}}.
\end{eqnarray}
Setting $l=K_{i+1}$, $t=\bar{\rho}_i$ and $s=\rho
_i$ for the inequality (\ref{3.8}) in \textbf{Step 1} and from (\ref{3.25}) we have

\begin{eqnarray}\label{3.26.0}
a_{i+1}&\leq&  S^{-\frac{p_-^\ast}{p_-}} 2^{p_-^\ast}    \{C_{\mu,K}[\int_{A_{K_{i+1},\rho_i}}(\frac{2^{i+3}}{R})^{p_-^\ast}|u-K_{i+1}|^{p_-^\ast}dx +(K_{i+1}^{p_+}+1)|A_{K_{i+1},\rho_i}|]\nonumber\\
&& + \int_{A_{K_{i+1}},\rho_i}(\frac{2^{i+1}}{R})^{p_-}|u-K_{i+1}|^{p_-}dx   \}^{\frac{p_-^\ast}{p_-}}\nonumber\\
&\leq &S^{-\frac{p_-^\ast}{p_-}} 2^{p_-^\ast}    \{C_{\mu,K}[\int_{A_{K_{i+1},\rho_i}}(\frac{2^{i+3}}{R})^{p_-^\ast}|u-K_{i+1}|^{p_-^\ast}dx +(K_{i+1}^{p_+}+1)|A_{K_{i+1},\rho_i}|] \nonumber\\
&&+ (\frac{2^{i+1}}{R})^{p_-}(\int_{A_{K_{i+1}},\rho_i}\frac{p_-}{p_-^\ast}|u-K_{i+1}|^{p_-^\ast}dx +\frac{p_-^\ast-p_-}{p_-^\ast} |A_{K_{i+1}},\rho_i|) \}^{\frac{p_-^\ast}{p_-}}\nonumber\\
&\leq &S^{-\frac{p_-^\ast}{p_-}} 2^{p_-^\ast}    \{C_{\mu,K}[(\frac{2^{i+3}}{R})^{p_-^\ast}\int_{A_{K_{i+1},\rho_i}}|u-K_{i+1}|^{p_-^\ast}dx +(K_{i+1}^{p_+}+1)|A_{K_{i+1},\rho_i}|] \nonumber\\
&&+ (\frac{2^{i+1}}{R})^{p_-}(\int_{A_{K_{i+1}},\rho_i}|u-K_{i+1}|^{p_-^\ast}dx + |A_{K_{i+1}},\rho_i|) \}^{\frac{p_-^\ast}{p_-}},\nonumber\\
\end{eqnarray}
where $C_{\mu,K}= (C'+C'' \mu K^{q_+-p_+})=(2C_{p,\varepsilon}+C\lambda 2^{p_+}+\mu K^{q_+-p_+}2^{p_++1}+1)$. Since $u>K_{i+1}>K_i...\geq K_0\geq 1$ on $A_{K_{i+1},\rho_i}$ and the definition of $K_i$, we have
\begin{eqnarray}\label{3.27}
(K_{i+1}-K_i)^{p_-^\ast}|A_{K_{i+1},\rho_i}|\leq\int_{A_{K_{i+1},\rho_i}}|u-K_{i+1}|^{p_-^\ast}dx\leq a_i,
\end{eqnarray}
\begin{eqnarray}\label{3.28}
|A_{K_{i+1},\rho_i}|\leq (\frac{2^{i+2}}{K})^{p_-^\ast}a_i\leq\frac{2^{(i+2)p_-^\ast}}{K_0^{p_-^\ast}}a_i\leq 2^{p_-^\ast(i+2)}a_i,
\end{eqnarray}
and
\begin{eqnarray}\label{3.29}
K_{i+1}^{p_+}|A_{K_{i+1},\rho_i}|\leq \frac{2^{p_-^\ast(i+2)}}{K_{i+1}^{p_-^\ast-p_+}}a_i\leq 2^{p_-^\ast(i+2)}a_i.
\end{eqnarray}
By these inequalities (\ref{3.26.0})--(\ref{3.29}), we have
\begin{eqnarray}\label{3.30}
a_{i+1}&\leq& S^{-\frac{p_-^\ast}{p_-}} 2^{p_-^\ast} \{C_{\mu,K}[(\frac{2^{i+3}}{R})^{p_-^\ast}a_i+2^{p_-^\ast(i+2)+1}a_i  ]+\frac{2^{(i+1)p_-}}{R^{p_-}} (a_i+2^{p_-^\ast(i+2)}a_i)\}^{\frac{p_-^\ast}{p_-}}\nonumber\\
&\leq&   S^{-\frac{p_-^\ast}{p_-}} 2^{p_-^\ast}[C_{\mu,K}(\frac{2^{3p_-^\ast}}{R^{p_-^\ast}}+2^{2p_-^\ast+1})+\frac{2^{p_-}}{R^{p_-}}+\frac{2^{p_-+2p_-^\ast}}{R^{p_-}}   ]^{\frac{p_-^\ast}{p_-}} 2^{\frac{2{p_-^\ast}^2i}{p_-}}   a_i^{\frac{p_-^\ast}{p_-}}\nonumber\\
&=& c (2^\frac{2{p_-^\ast}^2}{p_-})^i a_i^{1+\eta}\nonumber\\
&=&c b^i a_i^{1+\eta},
\end{eqnarray}
where $c= S^{-\frac{p_-^\ast}{p_-}} 2^{p_-^\ast}[C_{\mu,K}(\frac{2^{3p_-^\ast}}{R^{p_-^\ast}}+2^{2p_-^\ast+1}) +\frac{2^{p_-}}{R^{p_-}}+\frac{2^{p_-+2p_-^\ast}}{R^{p_-}}]^\frac{p_-^\ast}{p_-}$, $C_{\mu,K}=2C_{p,\varepsilon}+C\lambda 2^{p_+}+\mu K^{q_+-p_-}2^{p_+}+1$, $b=2^\frac{2{p_-^\ast}^2}{p_-}>1$ and $\eta=\frac{p_-}{N-p_-}$. From $u\in W^{1,p(x)}_0(\Omega)\hookrightarrow L^{p^\ast(x)}(\Omega)\subset L^{p_-^\ast}(\Omega)$, we get
\begin{equation}
\int_{A_{\frac{K}{2}},R}|u-\frac{K}{2}|^{p_-
^\ast}dx\leq \int_{A_{\frac{K}{2}},R}|u|^{p_-
^\ast}dx<+\infty,
\end{equation}
which implies $|u(x)|<+\infty$ a.e. in $A_{\frac{K}{2},R}$. That is $|A_{\frac{K}{2}},R|\rightarrow0$ as $K\rightarrow\infty$. Obviously,
\begin{equation}\label{3.32}
a_0=\int_{A_{\frac{K}{2}},R}|u-\frac{K}{2}|^{p_-
^\ast}dx\leq \int_{A_{\frac{K}{2}},R}|u|^{p_-
^\ast}dx\rightarrow 0,
\end{equation}
as $K\rightarrow\infty$. Next, we will find some suitable value of $K$ and $\mu$ such that the inequality
\begin{equation}\label{3.33}
a_0\leq c^{-\frac{1}{\eta}} b^{-\frac{1}{\eta^2}}
\end{equation}
holds in Lemma 2.4, where $b$, $c$ and $\eta$ are given by (\ref{3.30}), and $a_0=\int_{A_{\frac{K}{2}},R}|u-\frac{K}{2}|^{p_-
^\ast}dx$. This implies
\begin{eqnarray}
C_{\mu,K}( \frac{2^{3p_-^\ast}}{R^{p_-^\ast}}+2^{2p_-^\ast+1}   )+\frac{2^{p_-}}{R^{p_-}}+\frac{2^{p_-+2p_-^\ast}}{R^{p_-}}
\leq \frac{S 2^{-p_--2N}}{(\int_{A_{\frac{K}{2}},R}|u-\frac{K}{2}|^{p_-
^\ast}dx)^\frac{p_-}{N}}.\nonumber\\
\end{eqnarray}
From (\ref{3.32}) choose suitable $K$ to satisfy the following inequality
\begin{eqnarray}\label{3.26}
( \frac{2^{3p_-^\ast}}{R^{p_-^\ast}}+2^{2p_-^\ast+1}   )^{-1} [\frac{S 2^{-p_--2N}}{(\int_{A_{\frac{K}{2}},R}|u-\frac{K}{2}|^{p_-
^\ast}dx)^\frac{p_-}{N}}   -\frac{2^{p_-}}{R^{p_-}}-\frac{2^{p_-+2p_-^\ast}}{R^{p_-}}]-2C_{\varepsilon,p}-\lambda C2^{p_+}-1>0,
            \nonumber\\
\end{eqnarray}
and fix $\delta_1$ such that
\begin{eqnarray}
\mu \leq \delta_1&=&\frac{1}{K^{q_+-p_+}2^{p_++1}}\{( \frac{2^{3p_-^\ast}}{R^{p_-^\ast}}+2^{2p_-^\ast+1}   )^{-1} [\frac{S 2^{-p_--2N}}{(\int_{A_{\frac{K}{2}},R}|u-\frac{K}{2}|^{p_-
^\ast}dx)^\frac{p_-}{N}}   -\frac{2^{p_-}}{R^{p_-}}-\frac{2^{p_-+2p_-^\ast}}{R^{p_-}}]\nonumber\\
&&-2C_{\varepsilon,p}-\lambda C2^{p_+}-1\},
\end{eqnarray}
where $C_{p,\varepsilon}=\frac{p_+}{p_-}2^{p_+}(2p_+-2)^\frac{(p_+-1)p_-}{p_+}$. Thus we obtain $a_0\leq c^{-\frac{1}{\eta}}b^{-\frac{1}{\eta^2}}$ for $\mu\in [0,\delta_1]$. All the hypotheses of Lemma 2.4 are satisfied from (\ref{3.30}) and (\ref{3.33}), we obtain $a_i\rightarrow 0$ as $i\rightarrow\infty$, which implies
\begin{eqnarray}
u\leq K~\text{ on}~ B_{\frac{R}{2}},~\text{ for} ~\mu \in[0,\delta_1],
\end{eqnarray}
by (\ref{3.21.0}) and (\ref{3.22.0}). Set $\delta=\min\{\delta_0,\delta_1\}$, where $\delta_0>0$ is given by Lemma 2.5. Similarly, we can prove that $-u$ is also bounded on $B_{\frac{R}{2}}$. Since $u(x)|_{\partial\Omega}=0$ and $\Omega$ is a bounded smooth domain, we can obtain $|u|\leq K$ a.e. $\Omega$ for $\mu\in[0,\delta]$. We obtain the inequality (\ref{3.7}) and the proof is completed.

\end{document}